# Rationally connected varieties over local fields

By János Kollár*

## 1. Introduction

Let $X$ be a proper variety defined over a field $K$. Following [Ma], we say that two points $x, x' \in X(K)$ are *R-equivalent* if they can be connected by a chain of rational curves defined over $K$ (cf. (4.1)). In essence, two points are $R$-equivalent if they are "obviously" rationally equivalent. Several authors have proved finiteness results over local and global fields (cubic hypersurfaces [Ma], [SD], linear algebraic groups [CT-Sa], [Vo1], [Gi], [Vo2], rational surfaces [CT-Co], [CT1], [CT-Sk1], quadric bundles and intersections of two quadrics [CT-Sa-SD], [Pa-Su]).

$R$-equivalence is only interesting if there are plenty of rational curves on $X$, at least over $\bar{K}$. Such varieties have been studied in the series of papers [Ko-Mi-Mo1]–[Ko-Mi-Mo3]; see also [Ko1]. There are many *a priori* different ways of defining what "plenty" of rational curves should mean. Fortunately many of these turn out to be equivalent and this leads to the notion of *rationally connected* varieties. See [Ko-Mi-Mo2], [Ko1, IV.3], [Ko2, 4.1.2].

DEFINITION-THEOREM 1.1. *Let $\bar{K}$ be an algebraically closed field of characteristic zero. A smooth proper variety $X$ over $\bar{K}$ is called* rationally connected *if it satisfies any of the following equivalent properties*:

1. *There is an open subset $\emptyset \neq U \subset X$, such that for every $x_1, x_2 \in U$, there is a morphism $f : \mathbb{P}^1 \to X$ satisfying $x_1, x_2 \in f(\mathbb{P}^1)$.*

2. *For every $x_1, x_2 \in X$, there is a morphism $f : \mathbb{P}^1 \to X$ satisfying $x_1, x_2 \in f(\mathbb{P}^1)$.*

3. *For every $x_1, \ldots, x_n \in X$, there is a morphism $f : \mathbb{P}^1 \to X$ satisfying $x_1, \ldots, x_n \in f(\mathbb{P}^1)$.*

4. *Let $p_1, \ldots, p_n \in \mathbb{P}^1$ be distinct points and $m_1, \ldots, m_n$ natural numbers. For each $i$ let $f_i : \operatorname{Spec} \bar{K}[t]/(t^{m_i}) \to X$ be a morphism. Then there is a*

*Partial financial support was provided by the NSF under grant #DMS-9622394.

morphism $f : \mathbb{P}^1 \to X$ such that the Taylor series of $f$ at $p_i$ coincides with $f_i$ up to order $m_i$ for every $i$.

5. There is a morphism $f : \mathbb{P}^1 \to X$ such that $f^*T_X$ is ample.

6. For every $x_1, \ldots, x_n \in X$ there is a morphism $f : \mathbb{P}^1 \to X$ such that $f^*T_X$ is ample and $x_1, \ldots, x_n \in f(\mathbb{P}^1)$.

1.2.  The situation is somewhat more complicated in positive characteristic. The conditions of (1.1) are not mutually equivalent, but it turns out that (1.1.5) implies the rest ([Ko1, IV.3.9]). Such varieties are called *separably rationally connected*. The weakness of this notion is that there are even unirational varieties which are not separably rationally connected, and thus we do not cover all cases where finiteness is expected.

In characteristic zero, the class of rationally connected varieties is closed under smooth deformations ([Ko-Mi-Mo2]) and it contains all the known "rational like" varieties. For instance, unirational varieties and Fano varieties are rationally connected ([Na], [Ca], [Ko-Mi-Mo3]).

*Definition* 1.3.  By a *local field* $K$, we mean either $\mathbb{R}, \mathbb{C}, \mathbb{F}_q((t))$ or a finite extension of the $p$-adic field $\mathbb{Q}_p$. Each of these fields has a natural locally compact topology, and this induces a locally compact topology on the $K$-points of any algebraic variety over $K$, called the $K$-*topology*. The $K$-points of a proper variety are compact in the $K$-topology.

The aim of this paper is to study the $R$-equivalence classes on rationally connected varieties over local fields. The main result shows the existence of many rational curves defined over $K$. This in turn implies that there are only finitely many $R$-equivalence classes.

THEOREM 1.4.  *Let $K$ be a local field and $X$ a smooth proper variety over $K$ such that $X_{\bar{K}}$ is separably rationally connected. Then, for every $x \in X(K)$, there is a morphism $f_x : \mathbb{P}^1 \to X$ (defined over $K$) such that $f_x^*T_X$ is ample and $x \in f_x(\mathbb{P}^1(K))$.*

COROLLARY 1.5.  *Let $K$ be a local field and $X$ a smooth proper variety over $K$ such that $X_{\bar{K}}$ is separably rationally connected. Then*:

1. *Every $R$-equivalence class in $X(K)$ is open and closed in the $K$-topology.*

2. *There are only finitely many $R$-equivalence classes in $X(K)$.*

It is interesting to note that such a result should characterize rationally connected varieties.

CONJECTURE 1.6. *Let $X$ be a smooth proper variety defined over a local field $K$ of characteristic zero. Assume that $X(K) \neq \emptyset$ and there are only finitely many R-equivalence classes on $X(K)$. Then $X$ is rationally connected.*

*Note added in proof.* This conjecture was proved recently.

We show in (4.4) that (1.6) holds in dimensions 2 and 3. In general, it is implied by the geometric conjecture [Ko1, IV.5.6].

In the real case we can establish a precise relationship between the Euclidean topology of $X$ and the R-equivalence classes.

COROLLARY 1.7. *Let $X$ be a smooth proper variety over $\mathbb{R}$ such that $X_{\mathbb{C}}$ is rationally connected. Then the R-equivalence classes are precisely the connected components of $X(\mathbb{R})$.*

The following two consequences of (1.4) were pointed out to me by Colliot-Thélène. Yanchevskiĭ ([Ya]) proved that a conic bundle over $\mathbb{P}^1$ defined over a $p$-adic field $K$ is unirational if and only if it has $K$-points. (1.4) gives a new proof of this thanks to the fact that a conic is rational over a field $L$ if and only if it has an $L$-point. There are many other classes of varieties with a similar property, for each of which we obtain a unirationality criterion. (The precise general technical conditions are explained in (4.5).)

COROLLARY 1.8. *Let $K$ be a local field of characteristic zero and $X$ a smooth proper variety over $K$. Assume that there is a morphism $f : X \to \mathbb{P}^1$ whose geometric generic fiber $F$ is either*:

1. *a Del Pezzo surface of degree $\geq 2$,*

2. *a cubic hypersurface,*

3. *a complete intersection of two quadrics in $\mathbb{P}^n$ for $n \geq 4$, or*

4. *there is a connected linear algebraic group acting on $F$ with a dense orbit.*

*Then $X$ is unirational over $K$ if and only if $X(K) \neq \emptyset$.*

We also obtain a weaker result over global fields.

COROLLARY 1.9. *Let $\mathcal{O}$ be the ring of integers in a number field and $X$ a smooth proper variety defined over $\mathcal{O}$ satisfying one of the conditions (1.8.1–1.8.4). Assume that $X(\mathcal{O}) \neq \emptyset$. Then the mod $P$ reduction of $X$ is unirational over $\mathcal{O}/P$ for almost all prime ideals $P < \mathcal{O}$.*

*Remark* 1.10. The main theorem (1.4) and Corollary (1.8) hold for any field $K$ which has the property that on any variety with one smooth $K$-point the $K$-points are Zariski dense. Such fields are called *large fields* in [Po]. Examples of such fields are given in (2.3). If $K$ is real closed or complete with respect to a discrete valuation, then the $R$-equivalence classes are open and closed in the $K$-topology, but this does not imply the finiteness of $R$-equivalence classes unless $K$ is locally compact. In fact, finiteness fails in general, even for real closed fields.

*Acknowledgement.* I thank Colliot-Thélène for many helpful comments and references.

## 2. Smoothing lemmas

Let $f_t : \mathbb{P}^1 \to X$ be an algebraic family of morphisms for $t \neq 0$. As $t \to 0$, the limit of $f_t(\mathbb{P}^1)$ is a curve which may have several irreducible components. One can look at the limit as the image of a morphism from a reducible curve to $X$. (There does not seem to be a unique choice, though.) The *smoothing problem* studied in [Ko-Mi-Mo1]–[Ko-Mi-Mo3] attempts to do the converse of this. Given a reducible curve $C$ and a morphism $f_0 : C \to X$, we would like to write it as a limit of morphisms $f_t : \mathbb{P}^1 \to X$.

*Definition* 2.1. Let $C$ be a proper curve (possibly reducible and nonreduced) and $f : C \to X$ a morphism to a variety $X$. A *smoothing* of $f$ is a commutative diagram

$$\begin{array}{ccccc} C & \subset & S & \stackrel{F}{\to} & X \times T \\ \downarrow & & h \downarrow & & \downarrow \\ 0 & \in & T & = & T, \end{array}$$

where $0 \in T$ is a smooth pointed curve, $h : S \to T$ is flat and proper with smooth generic fiber, $C \cong h^{-1}(0)$ and $F|_C = f$. In this case we can think of $f : C \to X$ as the limit of the morphisms $F_t : h^{-1}(t) \to X$ as $t \to 0$.

Let $p_i \in C$ be points. We say that the above smoothing *fixes the $f(p_i)$* if there are sections $s_i : T \to S$ such that $s_i(0) = p_i$ and $F \circ s_i : T \to X \times T$ is the constant section $f(p_i)$ for every $i$.

Assume that $f : C \to X$ is defined over a field $K$. We say that $f$ is *smoothable over $K$* (resp. that $f$ is *smoothable over $K$ fixing the $f(p_i)$*) if there is a smoothing where everything is defined over $K$ (resp. which also fixes the $f(p_i)$).

The following smoothing result was established in [Ko-Mi-Mo2, 1.2]. First write $C$ as the special fiber of a surface $S \to T$ and then try to extend the morphism $f$ to $S$. This leads to the theory of Hom-schemes, discussed for instance in [Gro, 221]. In many applications we also would like to ensure that the $F_t$ pass through some points of $X$ and these points vary with $t$ in a prescribed manner. (This is the role of $Z$ in the next proposition.) Moreover, we may also require $X$ to vary with $t$, (thus $X = Y_0$ in the next result).

If $U \to V$ is a morphism, then $U_v$ denotes the fiber over $v \in V$.

PROPOSITION 2.2. *Let $T$ be a Noetherian scheme with a closed point $0 \in T$ and residue field $K$. Let $h : S \to T$ be a proper and flat morphism and $Z \subset S$ a closed subscheme such that $h : Z \to T$ is also flat. Let $g : Y \to T$ be a smooth morphism. Let $f : S_0 \to Y_0$ be a $K$-morphism and $p : Z \to Y$ a $T$-morphism such that $f|_{Z_0} = p|_{Z_0}$. Assume that*

1. $H^1(S_0, f^*(T_{Y_0}) \otimes I_{Z_0}) = 0$.

*Then there are*

2. *a scheme $T'$ with a $K$-point 0,*

3. *an étale morphism $(0 \in T') \to (0 \in T)$, and*

4. *a $T'$-morphism $F : S \times_T T' \to Y \times_T T'$,*

*such that $F|_{S_0} = f$ and $F|_{Z \times_T T'} = p \times_T T'$.*

*Note.* The statement of [Ko-Mi-Mo2, 1.2] is not exactly the above one. The assumption $H^1(S_0, f^*(T_{Y_0}) \otimes I_{Z_0}) = 0$ guarantees that the scheme $\text{Hom}(S, Y, p)$ defined in [Ko-Mi-Mo2, 1.1] is smooth over $T$ at $[f]$. Thus it has an étale section through $[f]$; this is our $T'$.

If $X$ is a variety over a field $K$, then we would like to find morphisms $\mathbb{P}^1 \to X$ defined over $K$, so we need to find $K$-points of $T'$. We do not have any control over $T'$ beyond those stated in (2.2). If $T$ is smooth over $K$ then $T'$ is also smooth and $0 \in T'$ is a $K$-point. This leads to the following question:

2.3. For which fields $K$ is it true that every curve with a smooth $K$-point contains a Zariski dense set of $K$-points? Characterizations of this property are given in [Po, 1.1].

The following are some interesting classes of such fields:

1. Fields complete with respect to a discrete valuation (This, in particular, includes the finite extensions of the $p$-adic fields $\mathbb{Q}_p$),

2. More generally, quotient fields of local Henselian domains,

3. $\mathbb{R}$ and all real closed fields,

4. Infinite algebraic extensions of finite fields and, more generally, pseudo algebraically closed fields (cf. [Fr-Ja, Chap. 10]).

In trying to apply (2.2), we see that the key point is to ensure the vanishing (2.2.1). In our cases the following easy lemma works.

LEMMA 2.4. *Let $C = \cup_{i=1}^m C_i$ be a proper curve whith only nodes. Set $C^i = \sum_{j=1}^i C_j$ and let $S_i := C_i \cap C^{i-1}$. Let $E$ be a vector bundle on $C$ such that $H^1(C_i, E|_{C_i} \otimes \mathcal{O}_{C_i}(-S_i)) = 0$ for every $i$. Then $H^1(C, E) = 0$.*

*Proof.* For every $i$ we have an exact sequence

$$0 \to E|_{C_i} \otimes \mathcal{O}_{C_i}(-S_i) \to E|_{C^i} \to E|_{C^{i-1}} \to 0.$$

By assumption $H^1(C_i, E|_{C_i} \otimes \mathcal{O}_{C_i}(-S_i)) = 0$; thus

$$H^1(C^i, E|_{C^i}) \cong H^1(C^{i-1}, E|_{C^{i-1}}),$$

and we are done by induction. □

2.5. It is easy to see that, given any proper nodal curve $C$, there is a smooth surface $h : S \to T$ whose central fiber is $C$. For us it will be easy and useful to construct $S \to T$ directly in each case.

## 3. Proof of the main theorem

*Definition* 3.1. A vector bundle $E$ over $\mathbb{P}^1$ is *ample* if $E \cong \sum_i \mathcal{O}(a_i)$ with $a_i > 0$ for every $i$. Equivalently, $E$ is ample if and only if $H^1(\mathbb{P}^1, E(-2)) = 0$. (Over $\mathbb{P}^1$ every vector bundle is a sum of line bundles; thus the equivalence of the two definitions is easy. Over other curves ampleness is defined very differently; see [Fu, p. 212].) By the upper semi continuity of cohomology groups, ampleness is an open condition.

More generally, let $E$ be a vector bundle over a surface $S$ and $g : S \to T$ a proper and flat morphism whose general fiber is $\mathbb{P}^1$. Let $D_1, D_2 \subset S$ be two Cartier divisors which are sections of $g$. Assume that

$$H^1(S_0, (E \otimes \mathcal{O}_S(-D_1 - D_2))|_{S_0}) = 0$$

for some $0 \in T$ where $S_0$ denotes the fiber over 0. ($S_0$ may be reducible and nonreduced.) Then $E|_{S_t}$ is ample for $0 \neq t \in T$ in an open neighborhood of 0.

First we prove a version of (1.4) which holds over any field.

THEOREM 3.2. *Let $K$ be a field and $X$ a smooth proper variety over $K$ such that $X_{\bar{K}}$ is separably rationally connected. Then for every $x \in X(K)$ there are*

1. *a smooth, affine, geometrically irreducible $K$-curve with a $K$ point $0 \in T'$, and*

2. *a $K$-morphism $\Phi : (T' \setminus \{0\}) \times \mathbb{P}^1 \to X$,*

*such that*

3. *$\Phi((T' \setminus \{0\}) \times \{(0 : 1)\}) = x$, and*

4. *$\Phi^* T_X|_{\{t\} \times \mathbb{P}^1}$ is ample for $t \neq 0$.*

*Proof.* Pick $x \in X(K)$. By (1.1.6) there is a morphism $g : \mathbb{P}^1_{\bar{K}} \to X_{\bar{K}}$ such that $g(0 : 1) = x$ and $g^* T_X$ is ample. $g$ is defined over a field extension $L = K(z) \supset K$ which we may assume to be Galois over $K$. (By composing $g$ with a suitable Frobenius one avoids inseparability problems.) Let $z_1 = z, z_2, \ldots, z_d$ denote the conjugates of $z$ and $g_i : \mathbb{P}^1_{\bar{K}} \to X_{\bar{K}}$ the corresponding conjugates of $g$ (with $g = g_1$). Note that $g_i(0 : 1) = x$ for every $i$.

Let $0 \in T$ be any smooth affine curve over $K$ with a $K$-point. Define the sections $s_0, \ldots, s_d : T \to T \times \mathbb{P}^1$ by $s_0(t) = (t, (0 : 1))$ and $s_i(t) = (t, z_i)$ for $i = 1, \ldots, d$. $s_0$ is defined over $K$ and the $s_i$ are defined over $L$. Let $S$ be the blow up of $T \times \mathbb{P}^1$ at the points $(0, z_1), \ldots, (0, z_d)$ and $h : S \to T \times \mathbb{P}^1 \to T$, the composite. Let $C_i$ denote the exceptional curve over the point $(0, z_i)$ and $C_0$ the birational transform of $0 \times \mathbb{P}^1$. Then $S_0 = h^{-1}(0) = C_0 + \cdots + C_d$ is a reduced curve. Its only singular points are the nodes $q_i := C_0 \cap C_i$. The sections $s_i$ lift to sections of $h$, these are denoted by $\bar{s}_i$. Let $Z \subset S$ be the image of the section $\bar{s}_0$. Set $p_i := \bar{s}_i(0) \in S_0$.

Observe that $S, h$ and $Z$ are defined over $K$ since the $z_i$ form a complete set of conjugates.

Fix a local parameter at $0 \in T$. This specifies local parameters at each $(0, z_i)$ and so gives $L$-isomorphisms $\tau_i : C_i \cong \mathbb{P}^1_L$ such that $\tau_i(q_i) = (0 : 1)$. The $\tau_i$ are conjugates of each other.

Define a morphism $f : S_0 \to X$ as follows. Set $f|_{C_i} = g_i \circ \tau_i$ for $i = 1, \ldots, d$ and let $f|_{C_0}$ be the constant morphism to $x$. These rules agree at the points $C_0 \cap C_i$, thus we get an $L$-morphism $f : S_0 \to X$. $f$ is in fact a $K$-morphism since the $g_i \circ \tau_i$ are conjugates of each other.

Set $Y = T \times X$ and define $p : Z \to Y$ by the rule $(t, \bar{s}_0(t)) \mapsto (t, x)$. $p$ is defined over $K$.

$f^* T_X|_{C_0} \cong \mathcal{O}_{C_0}^{\dim X}$, so $H^1(C_0, f^* T_X|_{C_0}(-[p_0])) = 0$. $f^* T_X|_{C_i}$ is ample for $i > 0$, so
$$H^1(C_i, f^* T_X|_{C_i}(-[p_i] - [q_i])) = 0.$$
Thus $H^1(S_0, f^* T_X(-\sum[p_i])) = 0$ by (2.4) and Proposition (2.2) applies.

Let $(0 \in T') \to (0 \in T)$ and $F : S \times_T T' \to X \times T'$ be as in (2.2). Let $0 \neq t \in T'(\bar{K})$ be a $\bar{K}$-point. The fiber of $S \times_T T'$ over $t$ is isomorphic to $\mathbb{P}^1_{\bar{K}}$ and the restriction of $F$ gives a $\bar{K}$-morphism $F_t : \mathbb{P}^1_{\bar{K}} \to X$ such that $F_t(0:1) = x$.

We have proved above that $H^1(S_0, f^*T_X(-\sum [p_i])) = 0$; thus
$$H^1(\mathbb{P}^1, F_t^*T_X(-d-1)) = 0$$
in a suitable Zariski open subset of $T'$ by (3.1). By shrinking $T'$, we may assume that it holds for every $t \in T' \setminus \{0\}$. So $F_t^*T_X$ is ample for $t \neq 0$.

Define $\Phi$ as the composite
$$\Phi : \mathbb{P}^1 \times (T' \setminus \{0\}) \xrightarrow{F} X \times T' \xrightarrow{\pi} X$$
where $\pi$ is the first projection. □

3.3. *Proof of* (1.4). Let $T'$ be as in (3.2). If $K$ is a large field, then we can pick a $K$-point $t \in T' \setminus \{0\}$. The induced morphism $f_x := \Phi_t$ has the required properties. □

## 4. Proof of the corollaries

*Definition* 4.1 ([Ma, §14]). Let $X$ be a variety over a field $K$. Two points $x, y \in X(K)$ are called *directly R-equivalent* if there is a $K$-morphism $f : \mathbb{P}^1 \to X$ such that $x = f(0:1)$ and $y = f(1:0)$. Two points $x, y \in X(K)$ are called *R-equivalent* if there is a sequence of points $x_0 = x, \ldots, x_m = y$ such that $x_i$ and $x_{i+1}$ are directly R-equivalent for $i = 0, \ldots, m-1$.

4.2. *Proof of* (1.5). If all the R-equivalence classes are open in the $K$-topology then they are also closed, since the complement of any equivalence class is the union of the other equivalence classes. $X(K)$ is compact in the $K$-topology, so there are only finitely many R-equivalence classes.

Let $U$ be an R-equivalence class. We need to prove that $U$ contains a $K$-open neighborhood for every $x \in U$. Let $f : \mathbb{P}^1_K \to X$ be a $K$-morphism such that $f(0:1) = x$. We may assume that $y := f(1:0) \neq x$.

Let $\mathrm{Hom}(\mathbb{P}^1, X, (1:0) \mapsto y)$ denote the universal family of those morphisms $f : \mathbb{P}^1 \to X$ such that $f(1:0) = y$ (cf. [Ko1, II.1.4]). Let $V \subset \mathrm{Hom}(\mathbb{P}^1, X, (1:0) \mapsto y)$ be an open subset containing $[f]$ such that $g^*T_X$ is ample for every $g \in V$. By [Ko1, II.3.5.3], this implies that the universal morphism
$$F : \mathbb{P}^1 \times \mathrm{Hom}(\mathbb{P}^1, X, (1:0) \mapsto y) \to X$$
is smooth away from $\{(1:0)\} \times \mathrm{Hom}(\mathbb{P}^1, X, (1:0) \mapsto y)$. The analytic inverse function theorem (cf. [Gra-Re, p. 102]) implies that a smooth morphism is open in the $K$-topology. Therefore the image
$$F(\mathbb{P}^1(K) \times \mathrm{Hom}(\mathbb{P}^1, X, (1:0) \mapsto y)(K))$$
contains an open neighborhood of $x$.

4.3. *Proof of* (1.7). $\mathbb{P}^1(\mathbb{R})$ is connected, hence every $R$-equivalence class is connected over $\mathbb{R}$. An open and closed subset of $X(\mathbb{R})$ is a union of connected components. □

4.4. *Remarks about* (1.6). Let $X$ be a smooth, proper variety over $K$. Assume that there is an open set $X^0 \subset X$ and a proper a morphism $g : X^0 \to Z$ such that $Z_{\bar{K}}$ is not covered by rational curves. Then there are countably many subvarieties $W_i \subsetneq Z$ such that every rational curve on $Z_{\bar{K}}$ is contained in $\cup_i W_i$ (This follows from the fact that there are only countably many families of rational curves on any variety, (cf. Ko1, II.2.11]).) Thus if $x \in X^0$ is a point such that $g(x) \notin \cup_i W_i$, then the $R$-equivalence class of $x$ is contained in $g^{-1}(g(x))$. If $X(K) \neq \emptyset$, then there are continuously many such $R$-equivalence classes.

Conjecture [Ko1, IV.5.6] asserts that, if $X$ is not rationally connected, then there is such a morphism $g : X^0 \to Z$. If $\dim X \leq 3$, then the conjecture is true by [Ko-Mi-Mo2, 3.2], and thus (1.6) holds in dimensions $\leq 3$.

4.5. *Proof of* (1.8). If $X$ is unirational, then clearly $X(K) \neq \emptyset$. To see the converse, first we establish that $X$ is separably rationally connected. This is well known; see, for instance, [Ko1, IV.6].

If $X(K) \neq \emptyset$, then by (1.4) there is a morphism $g : \mathbb{P}^1 \to X$ whose image is not contained in a fiber of $f$. Pulling back $f : X \to \mathbb{P}^1$ by $g$, we obtain a $K$-variety $f' : X' \to \mathbb{P}^1$ with geometric generic fiber $F'$. Moreover, $f'$ has a section over $K$. It is sufficient to prove that $X'$ is unirational over $K$, or that $F'$ is unirational over $K(t)$.

All three cases listed in Corollary 1.8 are varieties with the property that if they are defined over a field $L$ with a "sufficiently general" $L$-point then they are unirational over $L$. (In each case, sufficiently general means: outside an a priori given closed subset.) For Del Pezzo surfaces, see [Ma, IV.7.8] and for cubic hypersurfaces see [op cit., II.2.9]. Complete intersections of two quadrics are treated in [CT-Sa-SD, I, Prop.2.3]; for almost homogeneous spaces this is a result of Chevalley and Springer (see [Bo, 18.2] or [Ko1, IV.6.9]). Over a local field we have many choices for the rational curve $g : \mathbb{P}^1 \to X$, so there is no problem with the "sufficiently general" condition.

Most of the proof works in any characteristic but there are occasional inseparability problems, especially in the almost homogeneous case. □

4.6. *Proof of* (1.9). The key point is again to find rational curves $f_{\mathcal{O}/P}$ : $\mathbb{P}^1_{\mathcal{O}/P} \to X_{\mathcal{O}/P}$ for almost all $P$ which are not contained in a fiber. Following the proof of (3.3) we obtain $T'$ defined over $\mathcal{O}$. Let $T''$ be a smooth compactification of $T'$ and $B = \{0\} \cup (T'' \setminus T')$. Except for finitely many $P$, any $\mathcal{O}/P$-point of $T''_{\mathcal{O}/P} \setminus B_{\mathcal{O}/P}$ gives a desired rational curve.

By the Weil estimates (cf. [Ha, Ex. V.1.10]), a projective, geometrically irreducible curve of genus $g$ over $\mathbb{F}_q$ has at least $m$ points in $\mathbb{F}_q$ for $q - 2g\sqrt{q} + 1 \geq m$. Thus $T''_{\mathcal{O}/P} \setminus B_{\mathcal{O}/P}$ has points in $\mathcal{O}/P$ for almost all $P$ and the rest of the proof works as above. □

University of Utah, Salt Lake City, UT
*E-mail address*: kollar@math.utah.edu
*Current address*: Princeton University, Princeton, NJ
*E-mail address*: kollar@math.princeton.edu